\documentclass{agtart_a}
\pdfoutput=1

\usepackage{pinlabel}

\title{On links with cyclotomic Jones polynomials}

\author{Abhijit Champanerkar}
\givenname{Abhijit}
\surname{Champanerkar}
\address{Department of Mathematics and Statistics\\
University of South Alabama\\\newline
Mobile, AL 36688\\USA}
\email{achampanerkar@jaguar1.usouthal.edu}
\urladdr{}

\author{Ilya Kofman}
\givenname{Ilya}
\surname{Kofman}
\address{Department of Mathematics\\
College of Staten Island\\
City University of New York\\\newline
2800 Victory Boulevard\\Staten Island, NY 10314\\USA }
\email{ikofman@math.csi.cuny.edu}
\urladdr{}

\volumenumber{6}
\issuenumber{}
\publicationyear{2006}
\papernumber{58}
\startpage{1655}
\endpage{1668}

\doi{}
\MR{}
\Zbl{}


\received{5 June 2006}
\revised{}
\accepted{28 August 2006}
\published{14 October 2006}
\publishedonline{14 October 2006}
\proposed{}
\seconded{}
\corresponding{}
\editor{}
\version{}

\arxivreference{math.GT/0605631}

\AtBeginDocument{\let\bar\wbar\let\tilde\wtilde\let\hat\what}

\makeatletter
\def\cnewtheorem#1[#2]#3{\newtheorem{#1}{#3}[section]
\expandafter\let\csname c@#1\endcsname\c@thm}
\makeatother

\theoremstyle{plain}
 \newtheorem{thm}{Theorem}[section]
 \cnewtheorem{prop}[thm]{Proposition}
\newtheorem{lemma}[thm]{Lemma}
 \cnewtheorem{cor}[thm]{Corollary}
\theoremstyle{definition}
 \newtheorem{defn}{Definition}
 \makeautorefname{defn}{Definition}
 \newtheorem{example}{Example}

\renewcommand{\L}{\tilde{L}}
\newcommand{\Lx}{L_{\mbox{\psfig{file=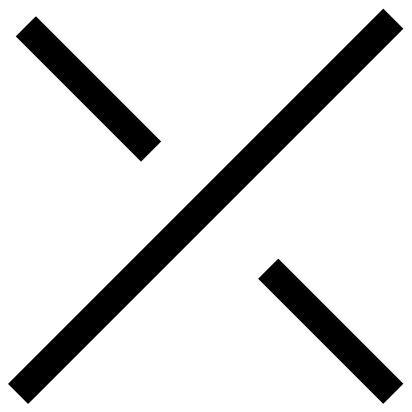,height=6pt}}}}

\renewcommand{\a}{\mathbf{a}}
\newcommand{\n}{\mathbf{n}}
\newcommand{\x}{\mathbf{x}}

\newcommand{\s}{\mathbf{s}}
\renewcommand{\r}{\mathbf{r}}

\newcommand{\kb}[1]{\ensuremath{\langle #1 \rangle}}

\newcommand{\e}{\varepsilon}
\newcommand{\myone}{\ensuremath{{\bf 1}}}

\renewcommand\theenumi{{\roman{enumi}}}

\newcommand{\asplice}{\ensuremath{\makebox[0.3cm][c]{\raisebox{-0.3ex}{\rotatebox{90}{$\asymp$}}}}}

\begin{document}

\begin{asciiabstract}
We show that if {L_n} is any infinite sequence of links with twist
number tau(L_n) and with cyclotomic Jones polynomials of increasing
span, then lim sup tau(L_n)=infty.  This implies that any infinite
sequence of prime alternating links with cyclotomic Jones polynomials
must have unbounded hyperbolic volume.  The main tool is the
multivariable twist--bracket polynomial, which generalizes the
Kauffman bracket to link diagrams with open twist sites.
\end{asciiabstract}

\begin{htmlabstract}
We show that if {L<sub>n</sub>} is any infinite sequence of links with
twist number &tau;(L<sub>n</sub>) and with cyclotomic Jones
polynomials of increasing span, then
lim&nbsp;sup&nbsp;&tau;(L<sub>n</sub>)=&infin;.  This implies that any infinite
sequence of prime alternating links with cyclotomic Jones polynomials
must have unbounded hyperbolic volume.  The main tool is the
multivariable <em>twist&ndash;bracket polynomial</em>, which generalizes
the Kauffman bracket to link diagrams with open twist sites.
\end{htmlabstract}

\begin{abstract}
We show that if $\{L_n\}$ is any infinite sequence of links with twist
number $\tau(L_n)$ and with cyclotomic Jones polynomials of increasing
span, then $\lim\sup\tau(L_n)=\infty$.  This implies that any infinite
sequence of prime alternating links with cyclotomic Jones polynomials
must have unbounded hyperbolic volume.  The main tool is the
multivariable {\em twist--bracket polynomial}, which generalizes the
Kauffman bracket to link diagrams with open twist sites.

\end{abstract}
\maketitle

\section{Introduction}
In \cite{ck05}, we showed that the Mahler measure of the Jones
polynomial converges under twisting for any link.  This is consistent
with the convergence of hyperbolic volume under the corresponding Dehn
surgery. 
In this note, we consider knots or links whose Jones polynomials are distinct but
with constant Mahler measure equal to one, which we call {\em cyclotomic}.
All known constructions of infinite sequences of hyperbolic knots or links
with the same Jones polynomial (eg,  Kanenobu \cite{Kanenobu1,Kanenobu2} and
Watson \cite{Watson})
require twisting at finitely many sites, and hence all such families have
bounded volume.
In contrast, any infinite sequence of hyperbolic alternating links
with cyclotomic Jones polynomials must have unbounded volume by
\fullref{coron} below.

For any link diagram $L$, two crossings are in the same {\em twist
class} if there is a simple closed curve that transversally intersects
the projection of $L$ only at the two crossing points and encloses
only adjacent bigons of the diagram.  The {\em twist number}
$\tau(L)$, originally defined by Lackenby in \cite{Lackenby}, is the
number of twist classes of crossings of $L$.
In general, for any collection of links with bounded twist number,
their Jones polynomials have bounded Mahler measure (Silver, Stoimenow and Williams
\cite{SSW}); see \fullref{Meik} below.
\begin{thm}\label{mainthm}
If $\{L(m)\}$ is any sequence of links with cyclotomic Jones
polynomials $V_{L(m)}(t)$, and 
$\displaystyle\lim_{m\to\infty}{\rm span}(V_{L(m)}(t))=\infty$, then
$\displaystyle\lim\sup\tau(L(m))=\infty$.
\end{thm}
Our theorem, together with Theorem 1 of \cite{Lackenby}, implies the
following corollary, where Vol$(S^3\setminus L)$ denotes hyperbolic
volume:
\begin{cor}\label{coron}
If $\{L(m)\}$ is any infinite sequence of distinct prime alternating links
with cyclotomic Jones polynomials, then $\lim\sup {\rm
Vol}(S^3\setminus L(m))=\infty$.
\end{cor}

A simple example that illustrates \fullref{mainthm} is the
sequence of connect sums of the figure--eight knot.  Prime knots with
cyclotomic Jones polynomials are relatively rare: there are 17 such
knots in 1.7 million knots with up to 16 crossings (see Remark 1 of
\cite{ck05}).  We do not know such an infinite sequence of hyperbolic
knots, although the simplest hyperbolic knots have Jones polynomials
with unusually small Mahler measure \cite{ckp}.

The main tool in the proof is our generalization of the Kauffman
bracket polynomial for links to the {\em twist--bracket polynomial}
$P(A,x_1,\ldots,x_k)$ for link diagrams with $k$ open twist sites.
After normalizing, we obtain a regular isotopy invariant for these
diagrams, which are sometimes called $2$--strand block diagrams.
Details and examples are given in the next section.

The twist--bracket polynomial is a special case of the multivariable
polynomial introduced in \cite{ck05}, which was defined using minimal
central idempotents in the Temperley--Lieb algebra.  By twisting on
only two strands at every twist site, we obtain an explicit formula in
terms of the Kauffman bracket.  In general, a link diagram with $k$
open twist sites gives rise to a $k$--linear form on the Temperley--Lieb
algebra, and the twist--bracket polynomial is naturally associated to
this form.  For details, see \cite{ck05,Watson}.

\subsubsection*{Acknowledgments}
The authors acknowledge support from the National Science Foundation under an FRG grant  
(DMS-0455978 and DMS-0456227, respectively).

\section{The twist--bracket polynomial}

For any link diagram $L$, its Kauffman bracket $\kb{L}\in\Z[A^{\pm 1}]$
equals the Jones polynomial $V_L(t)$, up to a monomial that depends on
the writhe $w(L)$ of the diagram:
If $t=A^{-4},\ V_L(t)=(-A)^{-3w(L)}\kb{L}$.

A {\em wiring diagram} is a regular $4$--valent planar graph with two
kinds of vertices, called twist sites, which are oriented either
horizontally or vertically.  We say that a link diagram $L$ can be
obtained from a wiring diagram $\L$ if $L$ is realized by inserting a
twist class of crossings at every twist site of $\L$, according to its
horizontal or vertical orientation.
For any $L$, there exists $\L$ with $v(\L)=\tau(L)$ such that $L$ can
be obtained from $\L$.
If all but $k$ twist sites have crossings inserted, we will say that
the diagram $\L$ has $k$ {\em open twist sites}.  A twist site will be
called {\em nugatory} if the isotopy class of the link does not depend
on how many crossings are inserted there.

Let $\L$ be a diagram with $k$ open twist sites.  For $\s \in
\{0,1\}^k$, let $\L_{\s}$ be obtained as follows: If the $i$--th twist
site is oriented vertically, insert $\asplice$ when $s_i=0$, and
$\asymp$ when $s_i=1$.  Otherwise, insert $\asymp$ when $s_i=0$, and
$\asplice$ when $s_i=1$.  For $\n \in \Z^k$, let $L_{\n}$ be the
diagram obtained from $\L$ by inserting at the $i$--th twist site:
$n_i$ half--twists on $2$ strands, which is $|n_i|$ crossings, with
sign according to right-- or left--handed twisting.  To be precise, if
$n_i=1$, we insert the crossing whose $A$--smoothing corresponds to
$s_i=0$.  If $n_i=0$, we insert $\asplice$ or $\asymp$ as in $s_i=0$.
\begin{center}\small\hair5pt
\labellist
\pinlabel $s_i=0$ [t] at 172 162
\pinlabel $s_i=1$ [t] at 298 162
\pinlabel $s_i=0$ [t] at 663 231
\pinlabel $s_i=1$ [t] at 860 231
\pinlabel $n_i=1$ [l] at 182 26
\pinlabel $n_i=1$ [l] at 652 26
\endlabellist
\includegraphics[width=4 in]{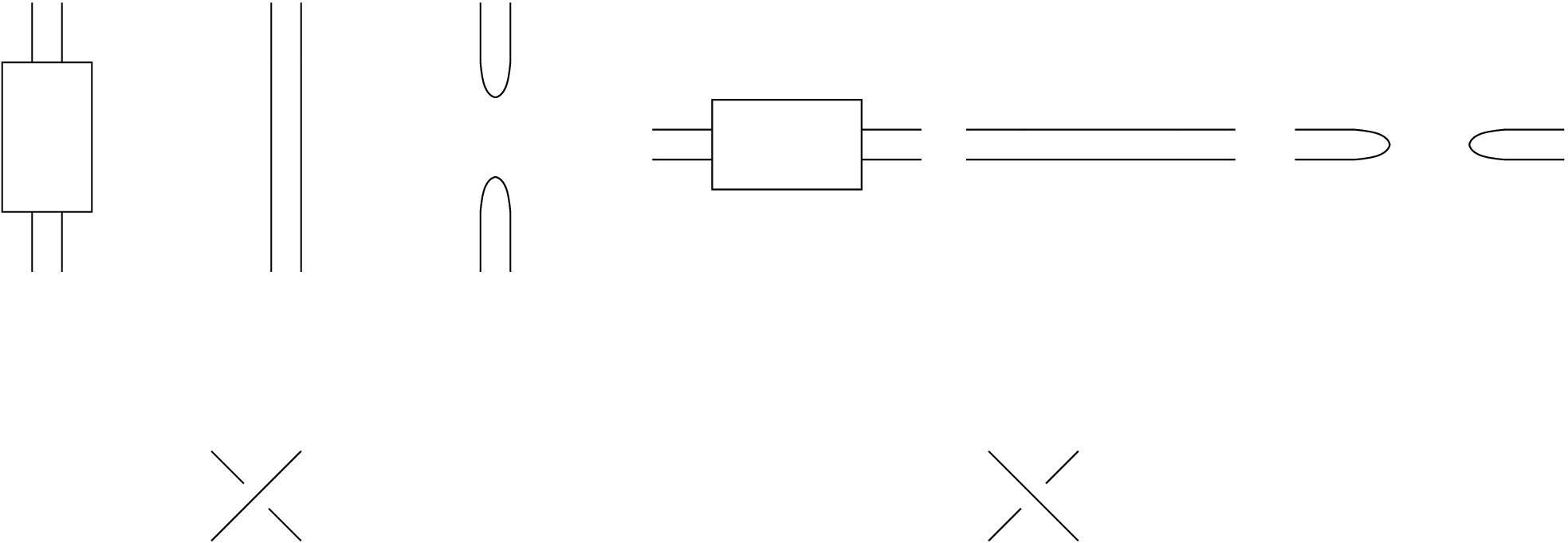}
\end{center}

Let $|\s|=\sum_{i=1}^{k} s_i$ and $\sigma(\n)=\sum_{i=1}^{k} n_i$.
Let $\delta=-A^2-A^{-2}$.

\begin{defn}\label{PAxx}
For any $\L$ with $k$ open twist sites, we define its {\em twist--bracket} by:
$$ P(A,x_1,\ldots,x_k) = \sum_{\s\in\{0,1\}^k} \left( \prod_{i=1}^k (x_i - 1)^{s_i}\;\delta^{1-s_i}\;\right)\kb{\L_{\s}} $$
\end{defn}
From the definition, $P(A,x_1,\ldots,x_k)\in \Z[A^{\pm 1},x_1,\ldots,x_k]$.

\begin{prop}\label{PAnn}
If $\L$ is any diagram with $k$ open twist sites, then
$$\delta^k\, \kb{L_{\n}}=A^{\sigma(\n)}\ P(A,(-A^{-4})^{n_{1}},\ldots,(-A^{-4})^{n_{k}}). $$
\end{prop}
\begin{proof}
First, suppose $\L$ has only one open twist site, which we orient vertically.
We claim that, similar to Proposition 3.3 \cite{ck05},
\begin{equation}\label{twtwist}
\delta \,\kb{L_n}= A^{n}\Big(\delta\,\kb{L_{\asplice}}+((-A^{-4})^n-1)\,\kb{L_\asymp}\Big).
\end{equation}
If $\{\myone, e_1\}$ is the $TL_2$ basis,
$p_0=\myone-e_1/\delta$ is the Jones--Wenzl idempotent, and
$p_1=\myone - p_0$ is the orthogonal idempotent.
Let $\rho:\thinspace B_2\rightarrow TL_2$ be given by $\rho(\sigma_1)=A\, \myone+A^{-1}\, e_1$.
Following the proof of Proposition 3.3 \cite{ck05}, with $\Delta=\sigma_1$,
$\rho(\Delta)=A p_0 -A^{-3}p_1$, hence $\rho(\Delta^{n})=
A^n p_0+(-1)^n A^{-3n}p_1$.
\begin{eqnarray*}
\rho(\Delta^{n})
&=& A^n (p_0+(-1)^n A^{-4n}p_1) \\
&=&A^{n}(\myone - \frac{e_1}{\delta} +(-1)^n A^{-4n}\frac{e_1}{\delta})\\
&=& A^{n}\Big(\myone + \frac{((-1)^nA^{-4n}-1)}{\delta}\,e_1) \Big)
\end{eqnarray*}
Using the bilinear form on $TL_2$, $\kb{L_n}=\langle \L,\Delta^n \rangle$, so \eqref{twtwist} follows.

If we repeatedly apply \eqref{twtwist} at each twist site, we obtain:
$$ \delta^k\,\kb{L_{\n}}=A^{\sigma(\n)}\sum_{\s\in\{0,1\}^k}
\left(\prod_{i=1}^k\left((-A^{-4})^{n_i}-1\right)^{s_i}\;\delta^{1-s_i}\;\right) \kb{\L_{\s}}\proved $$
\end{proof}

\begin{example}\label{exon}\rm
Double twist links $L(m,n)$ are obtained from a wiring diagram $\L$
with two twist sites such that $\kb{\L_{(0,0)}}=\kb{\L_{(1,1)}}=1$ and
$\kb{\L_{(0,1)}}=\kb{\L_{(1,0)}}=\delta$. We get
$P(A,x,y)=\delta^{2}(x+y-1)+(x-1)(y-1)$.
\end{example}

\begin{example}\rm
2--bridge links $L(n_1,\ldots,n_k)$ are obtained from a wiring diagram $\L^k$ 
that depends on the parity of $k$ (see eg, Lickorish \cite[Figure 1.8]{Lickorish}).
$\L^k$ satisfies the following recurrence:
$\L^k_{(\s,1)}=\L^{k-1}_{\s}$,\, $\L^k_{(\s,0,0)}=\L^{k-2}_{\s}$,\, $\L^k_{(\s,1,0)}=\bigcirc
\cup \L^{k-2}_{\s}$.
If $P_k(A,x_1,\ldots,x_k)$ is the twist--bracket for $\L^k$, then
$P_k$ satisfies the following recurrence:
$$P_k=(x_k-1)\, P_{k-1}\, +\, \delta^2 \, x_{k-1}P_{k-2}$$
with $P_1=\delta^2 + (x_1-1)$, which is the twist--bracket for (2,n)
torus links, and $P_2=\delta^{2}(x_1+x_2-1)+(x_1-1)(x_2-1)$, which is
the same twist--bracket as in \fullref{exon}.
\end{example}

\begin{example}\rm
Pretzel links $\mathcal{P}(n_1,\ldots,n_k)$ are obtained from a wiring
diagram $\L$ with $k$ twist sites such that
$\kb{\L_{\s}}=\delta^{k-|\s|-1}$ for $\s\neq(1,\ldots, 1)$ and
$\kb{\L_{(1,\ldots,1)}}=\delta$.
The following gives a closed formula for the bracket of any pretzel link:
$$P(A,x_1,\dots,x_k)=\frac{1}{\delta}\left(\prod_{i=1}^k\;(x_i-1+\delta^2) +
(\delta^2-1)\prod_{i=1}^k\;(x_i-1)\right)$$
\end{example}

\begin{example}\rm
Kanenobu links $K(p,q)$ given in \cite{Kanenobu1} are obtained from a
diagram $\L$ with two open twist sites such that
$\kb{\L_{(1,0)}}=\kb{\L_{(0,1)}}=\delta$, $\kb{\L_{(1,1)}}=\delta^2$ and
$\kb{\L_{(0,0)}}=\kb{8_9}=\kb{4_1}^2$, where both knot diagrams have
writhe $=0$.  After simplifying, we get
$P(A,x,y)=\delta^2(\kb{8_9}-1+xy)$.  For $p+q=$ constant, the writhe is
constant, so these links have the same Jones polynomial.
\end{example}

\begin{example}\rm
Kanenobu links $K(n_1,\ldots,n_k)$ given in Figure 9 of \cite{Kanenobu2} are
obtained from a diagram $\L$ with $k$ open twist sites such
that $\kb{\L_{\s}}=\delta^{|\s|}$ for $\s\neq (0,\ldots,0)$. After simplifying, we get:
$$P(A,x_1,\ldots,x_k)=\delta^k\left(\kb{\L_{(0,\ldots,0)}}-1+\prod_{i=1}^k x_i\right)$$ 
For $\sum n_i=$ constant, the writhe is constant, so these links
have the same Jones polynomial.
\end{example}

If $\L$ has only one open twist site (oriented vertically), \fullref{PAnn} implies
\begin{equation}\label{PAx}
\delta\,\kb{L_n}=A^n P(A,(-A^{-4})^n)
\quad {\rm with } \quad
P(A,x) =\delta\,\kb{L_{\asplice}}+ (x-1) \kb{L_\asymp}.
\end{equation}

\begin{prop}\label{xprod}
If $P(A,x) = x\cdot f(A)$ for some $f$, then $V_{L_n}(t) = V_{L_{\asplice}}(t)$
for all $n$.
\end{prop}

\begin{proof} By \eqref{PAx}, $P(A,x) = x\cdot f(A)$ occurs only if
$\delta\,\kb{L_{\asplice}}=\kb{L_\asymp}$.  This has several
implications:  First, by \eqref{twtwist},
$$\delta\,\kb{L_n} = A^{n}(-A^{-4})^n \kb{L_\asymp} =
(-A^{-3})^n\kb{L_\asymp}.$$
Second, at $A=1$, $\kb{L_\asymp}=(-2)\kb{L_{\asplice}}$.
Generally, if $\mu(L)$ is the number of components
of $L$, and $c(L)$ is the number of crossings of $L$, then at $A=1$,
\begin{equation}\label{Aon}
\kb{L}=(-1)^{w(L)}(-2)^{\mu(L)-1}=(-1)^{c(L)}(-2)^{\mu(L)-1}
\end{equation}
Therefore, the number of components changes: $\mu(L_\asymp)=\mu(L_{\asplice})+1$.
The following are possibilities for $L_{\asplice}$ at the twist site, up to
crossing changes in the rest of the link diagram:
$$\mbox{\psfig{file=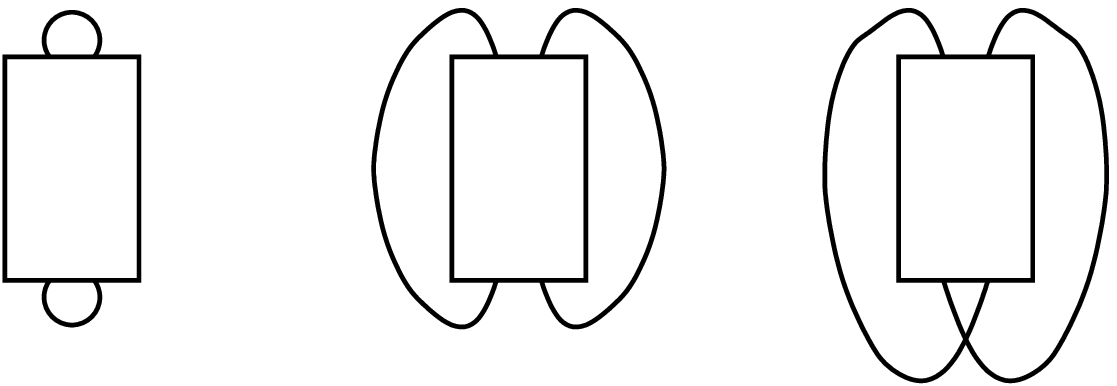,height=0.5in}}$$
The condition $\mu(L_\asymp)=\mu(L_{\asplice})+1$ excludes all except the first
type.  At such a twist site, the strands of $L_{\asplice}$ must be oriented in
opposite directions:
$$\mbox{\psfig{file=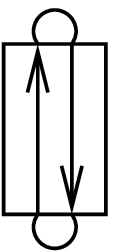,height=0.5in}}$$
Therefore, $w(L_n) = w(L_{\asplice}) - n$.
$$\eqalignbot{
V_{L_n}(t) &= (-A^{-3})^{w(L_n)}\kb{L_n} = (-A^{-3})^{w(L_{\asplice})-n} (-A^{-3})^n\kb{L_\asymp}/\delta \cr
&= (-A^{-3})^{w(L_{\asplice})}\kb{L_{\asplice}} = V_{L_{\asplice}}(t)}\proved
$$
\end{proof}

For any $\L$ with $k$ open twist sites, we define its {\em normalized twist--bracket} by
$$ \bar{P}_{\L}(A,x_1,\ldots,x_k) = \frac{1}{\delta^k}\,P_{\L}(A,x_1,\ldots,x_k).$$
The Kauffman bracket is an invariant of regular isotopy for link
diagrams.  We now show that the normalized twist--bracket is an invariant of
regular isotopy for link diagrams with open twist sites.
Two link diagrams with open twist sites will be called {\em regularly isotopic} 
if they are related by the following moves:
\begin{enumerate}
\item Planar isotopy.
\item Reidemeister moves $II$ and $III$, naturally augmented so that strands of the
  link can pass across open twist sites.
\item Combining and separating adjacent open twist sites:
$$\mbox{\psfig{file=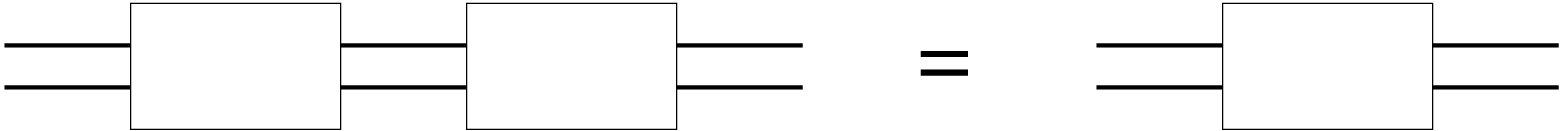,width=3in}}$$
\end{enumerate}

\begin{prop}
  If two link diagrams with open twist sites are regularly isotopic,
  then their normalized twist--brackets are equal, after possibly
  changing variables.
\end{prop}
\begin{proof} 
The invariance of the Kauffman bracket under moves (i) and (ii)
implies invariance of the twist--bracket under these moves.  We now
consider move (iii) for diagrams $\L^1$ and $\L^2$, which are
otherwise the same.  The terms of $\bar{P}_{\L^1}$ and $\bar{P}_{\L^2}$ for all
$\s$ in the twist sites shown are given in the following table:

\begin{tabular}{c|c|ccc}
 & 0 & 1 & & \\ 
\mbox{\psfig{file=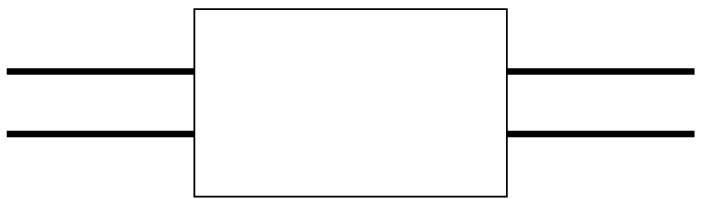,height=0.1in}} & $\kb{\L^1_0}$ & $\frac{(x-1)}{\delta}\kb{\L^1_1}$ \\ 
 &  &  & & \\ 
 & 00 & 10 & 01 & 11 \\ 
\mbox{\psfig{file=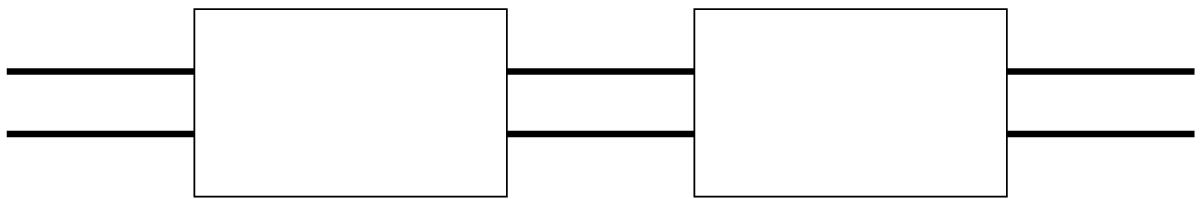,height=0.1in}} & $\kb{\L^2_{00}}$ & $\frac{(x_1-1)}{\delta}\kb{\L^2_{10}}$ & 
$\frac{(x_2-1)}{\delta}\kb{\L^2_{01}}$ & $\frac{(x_1-1)(x_2-1)}{\delta^{2}}\kb{\L^2_{11}}$ \\
\end{tabular}

We now observe that 
$\displaystyle \ \kb{\L^1_0}=\kb{\L^2_{00}}\quad {\rm and }\quad \kb{\L^1_1}=\kb{\L^2_{10}}=\kb{\L^2_{01}}=\frac{\kb{\L^2_{11}}}{\delta}$.

Since $(x_1-1)+(x_2-1)+(x_1-1)(x_2-1)=(x_1 x_2 - 1)$, we have that 
$$ \bar{P}_{\L^1}(A, (x_1 x_2)) = \bar{P}_{\L^2}(A, x_1, x_2) = \bar{P}_{\L^2}(A, (x_1 x_2), 1) = \bar{P}_{\L^2}(A, 1, (x_1 x_2)). $$ 
Hence, after possible variable changes, the 
normalized twist--bracket $\bar{P}_{\L}$ is an
invariant of $\L$ up to regular isotopy.  \end{proof}

\begin{lemma}  \label{sixthroot}
Let $\L$ be any diagram with $k$ open twist sites, obtained as above
from a wiring diagram with $N$ twist sites by inserting $n_j$
crossings at the $j$--th twist site for each $k<j\leq N$.
If $\xi=e^{2\pi i/6}$ then:
$$P(\xi,x_1,\ldots,x_k)=\prod_{i=1}^kx_i\prod_{j=k+1}^N(-\xi^{-3})^{n_j}$$
\end{lemma}
\begin{proof}
Since $\xi^4+\xi^2+1=0$, if $A=\xi$ then $\delta=-A^2-A^{-2}=1$.
If $\L'$ is a wiring diagram, then $\L_{\s}'$ is a collection of loops, so $\kb{\L_{\s}'}=\delta^{\mu-1}=1$.
Therefore: 
$$P(\xi,x_1,\ldots,x_N)=\sum_{\s\in\{0,1\}^N}\left(\prod_{i=1}^N(x_i - 1)^{s_i}\right)=\prod_{i=1}^N x_i$$
The last equality follows from the identity,
$$\displaystyle{\prod_{i=1}^N (X_i+1) =\sum_{\s\in\{0,1\}^N} \
\prod_{i=1}^N X_i^{s_i}}$$
with $X_i=x_i-1$. As $\L$ is obtained from $\L'$ by inserting $n_j$ crossings at the $j$--th twist
site for each $k<j\leq N$, the result follows by \fullref{PAnn}.
\end{proof}

\section{Mahler measure preliminaries}
Let $f \in \C[z^{\pm 1}_1,\ldots,z^{\pm 1}_s]$.  The Mahler measure of $f$ is defined by
$$ M(f) = \exp\int_0^1\cdots\int_0^1\log\left|f(e^{2\pi i\theta_1},\ldots,e^{2\pi i\theta_s})\right|
{\rm d} \theta_1\cdots {\rm d} \theta_s. $$
If $s=1,\; f(z)=a_0 z^k\prod_{i=1}^n(z-\alpha_i)$, then by Jensen's formula,
$$ M(f) = |a_0|\prod_{i=1}^n {\rm max}(1,|\alpha_i|). $$
Thus, for any $k\in \Z\setminus 0$, $M(f(z))=M(f(z^k))$.  Also, from the definition we get
\begin{equation}\label{Mf}
M(f(z_1,\ldots,z_i,\ldots,z_s))=M(f(z_1,\ldots,\pm z_i,\ldots,z_s)).
\end{equation}   
The Mahler measure is multiplicative, $M(f_1 f_2) = M(f_1) M(f_2)$,
so it can be naturally extended to rational functions of
Laurent polynomials.

A monic irreducible polynomial with coefficients in $\Z$ is called
{\em cyclotomic} if all of its zeros are primitive roots of
unity. Hence, the Mahler measure of a cyclotomic polynomial equals
one.  Multivariable polynomials whose Mahler measure equals one can be
expressed in terms of cyclotomic polynomials as follows:

\begin{lemma}[Boyd \cite{Boyd81}]\label{cpol} If $F\in \Z[x_1^{\pm1},\ldots, x_s^{\pm1}]$, then
$M(F)=1$ if and only if there are cyclotomic polynomials $\phi_j$ such that:
$$F=\prod_{i=1}^{s}x_i^{a_{i}}\prod_{j=1}^{m}\phi_j\left(\prod_{i=1}^{s}x_i^{b_{i,j}}\right)$$
\end{lemma}
Henceforth, we will use {\em cyclotomic} more loosely to mean constant Mahler measure equal to one. 

For a vector $\x\in\Z^s$, let $h(\x)=\max|x_i|$ and
$$\nu(\x) = \min\{ h(\a)\;|\; \a\in\Z^s\backslash \{ \mathbf{0}\},\; \a\cdot\x=0 \}. $$
For example, $\nu(1,d,\ldots,d^{s-1})=d$.
\begin{lemma}[Boyd, Lawton \cite{Lawton}]\label{BL}
For every $f \in \C[z^{\pm 1}_1,\ldots,z^{\pm 1}_s]$,
\[ M(f) = \lim_{\nu(\x)\to\infty} M\left(f(z^{x_1},\ldots,z^{x_s})\right) \]
with the following useful special case: $M(f) = \lim\limits_{d\to\infty}
M\left(f(z,z^d,\ldots,z^{d^{s-1}})\right)$.
\end{lemma}
We will use the following slight extension of \fullref{BL}: 
\begin{lemma}\label{BLe}
For every $f \in \C[z^{\pm 1}_1,\ldots,z^{\pm 1}_s]$ and every $\e\in\{\pm 1\}^s$,
\[ M(f) = \lim_{\nu(\x)\to\infty} M\left(f((\e_1 z)^{x_1},\ldots,(\e_s z)^{x_s})\right). \]
\end{lemma}
\begin{proof} 
Fix $f$.  Let $\displaystyle \psi(\x,\e) = M\left(f(\e_1 z^{x_1},\ldots,\e_s z^{x_s})\right)$.
By \eqref{Mf} and \fullref{BL}, for any $\e\in\{\pm 1\}^s$,
$$ M(f) = M(f(\e_1 z_1,\ldots,\e_s z_s)) = \lim_{\nu(\x)\to\infty}\psi(\x,\e).$$
Since there are finitely many $\e$, this implies  
$$ \lim_{\nu(\x)\to\infty} \min_{\e} \psi(\x,\e) = M(f)=\lim_{\nu(\x)\to\infty} \max_{\e} \psi(\x,\e). $$
The result now follows by the Squeeze Theorem.
\end{proof}

\begin{lemma} \label{nu}
Let $\{ \x_m \}$ be an infinite sequence of vectors in $\Z^k$. If
$\nu(\x_m) \leq R$ for all $m$, then there exists an infinite
subsequence $\{m_j\}$ and $\a \in \Z^k \backslash \{ \mathbf{0}\}$
such that $\x_{m_j}\cdot \a =0$ for all $j$.
\end{lemma}
\begin{proof} Since $\nu(\x_m)$ is a bounded nonnegative integer for all $m$, we
can pass to a subsequence $\{m_i\}$ such that $\nu(\x_{m_i})=C$ for
all $i$. Hence for each $i$, there exists $\a_{m_i} \in \Z^k
\backslash \{ \mathbf{0}\}$ such that $\a_{m_i}\cdot \x_{m_i}=0$ and
$h(\a_{m_i})=C$. Thus for all $i,\ \a_{m_i} \in [-C,C]^k \cap
\Z^k \backslash \{\mathbf{0}\}$. Since there are only
finitely many integer lattice points in $[-C,C]^k$, there exists a subsequence
$\{\a_{m_j} \}$ that is constant. \end{proof}

If $\{ \x_m \}$ is a sequence of vectors with bounded $\nu(\x_m)$,
\fullref{nu} allows us to pass to a subsequence and apply
\fullref{BL}, though the limiting polynomial will have fewer
variables.

\begin{example}\rm Let $\x_m=(1,m,2m+1)$.  Since $(1,2,-1)\cdot
(1,m,2m+1)=0,\ \nu(\x_m) \leq 2$.  Let $P(t,x,y)$ be a $3$--variable
polynomial.  As $\nu(\x_m) \leq 2$, it may not be true that
$M(P(t,t^m,t^{2m+1}))\to M(P(t,x,y))$ as $m\to \infty$. Let
$Q(t,x)=P(t,x,tx^2)$, then $Q(t,t^m)=P(t,t^m,t^{2m+1})$ and $\nu(1,m)
\to \infty$ as $m \to \infty$. Now by \fullref{BL}, we can conclude
that $M(P(t,t^m,t^{2m+1}))\to M(Q(t,x))$ as $m\to \infty$. The linear
dependence of the components $(1,m,2m+1)$ gives us a limit in one
fewer variable.
\end{example}

For $F\in \Z[x_1^{\pm1},\ldots, x_s^{\pm1}]$, $M(F) \leq ||F||$ where $||F||$ denotes
 the $L^2$ norm of coefficients of $F$ (see Schinzel \cite{Schinzel} for more details).

\begin{prop}\label{Meik} 
Let $L$ be a link diagram with $\tau(L)=k$ then $M(V_L(t))\leq 8^k$.
\end{prop}
\begin{proof} As $V_L(A^{-4})=(-A)^{-3w(L)}\kb{L}$, $M(V_L(t))=M(\kb{L})$. 
Since $M(\delta)=1$, by \fullref{PAnn}, $M(V_L(t))=M(\kb{L})=M(P(A,(-A^{-4})^{n_1},\ldots, (-A^{-4})^{n_k}))$.
\begin{eqnarray*}
M(V_L(t)) &\leq& ||P(A,(-A^{-4})^{n_1},\ldots, (-A^{-4})^{n_k})|| \leq ||P(A,x_1,\ldots,x_k)||\\
&\leq& \sum_{\s\in\{0,1\}^k} \left( \prod_{i=1}^k ||(x_i - 1)||^{s_i}\;||\delta||^{1-s_i}\;\right) ||\kb{\L_{\s}}||\\
&\leq& \sum_{\s\in\{0,1\}^k} \left( \prod_{i=1}^k 2^{s_i}\;2^{1-s_i}\;\right)2^k = 8^k \\
\end{eqnarray*}
The last inequality follows from the fact that $\kb{\L_{\s}}=\delta^r$ where $0\leq r \leq k$. \end{proof}

\section{Cyclotomic Jones polynomials}

Let $L$ be any link diagram.  Since $V_L(t)=(-A)^{-3w(L)}\kb{L}$ with
$t=A^{-4}$, the Mahler measure is unchanged: $M(V_L(t))=M(\kb{L})$.
The span of a Laurent polynomial is the difference of its highest and
lowest degrees, so $4\, {\rm span} (V_L(t))={\rm span}(\kb{L})$.
Let $L_m$ be the diagram obtained from $L$ by inserting $m$
half--twists on $2$ strands of $L$.
\begin{prop}\label{Propon}
Distinct $V_{L_m}(t)$ are cyclotomic for only finitely many $m$.
\end{prop}
\begin{proof}
We consider $L$ as having one open twist site, with the original link given by $L_{\asplice}$.
By \eqref{PAx}, $\delta\, \kb{L_m}=A^m P(A,(-A^{-4})^m)$, so then
$$ M(V_{L_{m}}(t))= M(\kb{L_m}) = M(P(A,(-A^{-4})^m)). $$ If
$M(V_{L_m}(t))=1$ for infinitely many $m$ then, after passing to a
subsequence, $M(P(A,x))= \lim_{m\to\infty} M(P(A, A^{-4m})) =1$ by
\fullref{BLe}.  Hence by \fullref{cpol}, $P(A,x)$ equals up to
monomials the product of cyclotomic polynomials evaluated at monomials
in $A$ and $x$.  Since the Jones polynomials vary with $m$, by
\fullref{xprod}, $P(A,x) \neq x\cdot f(A)$ for any $f$.  By
\eqref{PAx}, $P(A,x)$ is linear in $x$, so there is exactly one
cyclotomic linear factor evaluated at $A^r x$. Hence,
$$P(A,x)=(1\pm A^rx)A^\alpha\prod_{j=1}^{\ell}\phi_j(A^{m_j})\quad \Rightarrow\quad
\frac{\partial P}{\partial x}(A,x) = \pm A^r\ P(A,0). $$
Therefore by \eqref{PAx},
$ \kb{L_\asymp} = \pm A^r \left(\delta\,\kb{L_{\asplice}}-\kb{L_\asymp}\right) $.

At $A=1$, the equation above reduces to
$$ \kb{L_\asymp}\ = \pm\left( -2\kb{L_{\asplice}} - \kb{L_\asymp}\right). $$ If
the sign is negative, then $\kb{L_{\asplice}}=0$, and if the sign is positive,
then $\kb{\Lx}=\kb{L_{\asplice}}+\kb{L_\asymp}=0$, both of which are
contradictions by \eqref{Aon}.  \end{proof}

\begin{proof}[Proof of \fullref{mainthm}]
Let $\{ L(m)\}$ be an infinite sequence of link diagrams whose Jones
polynomials have increasing span and are cyclotomic: For all $m$,
$M(\kb{L(m)})=1$.  Suppose that $\{L(m)\}$ has bounded twist
number. Passing to a subsequence if necessary, we can assume that all
the links are obtained from the same wiring diagram with $N$ twist
sites, and since the Jones polynomials change, that these are
non-nugatory twist sites.  Hence, there is a sequence
$\{\n(m)=(n_1(m),\ldots,n_N(m) )\}$ in $\Z^N$ such that
$L(m)=L_{\n(m)}$.

Let $|\n(m)|=\sum_{i=1}^{N}|n_i(m)|$, which is the number of crossings
of $L_{\n(m)}$.  For any $L$, span$(V_L(t))$ is a lower bound for the
crossing number of $L$, so $|\n(m)|\to\infty$.  Passing to a
subsequence if necessary, we can assume that for all $i$, either
$n_i(m)\to\infty$, or $n_i(m)$ is constant for all $m$.  Hence, this
subsequence $\{ L(m)\}$ is obtained from a fixed diagram $\L$ with
$k\leq N$ open non-nugatory twist sites at which $n_i(m)\to\infty$.
We will call these {\em active twist sites}.  By \fullref{Propon}, to have an infinite sequence of links with distinct
cyclotomic Jones polynomials, there must be at least two active twist
sites.

Suppose $\{ L(m)\}$ has active twist sites for $1\leq i\leq k$, and
$P(A,x_1,\ldots,x_k)$ is given by \fullref{PAxx}.  By
\fullref{PAnn}, $V_{L_{\n(m)}}(t)$ is determined by setting
$x_i=(-A^{-4})^{n_i}$.  Since $4\, {\rm span} (V_L(t))={\rm
span}(\kb{L})$, our hypothesis that ${\rm
span}(V_{L_{\n(m)}}(t))\to\infty$ implies
\begin{equation}\label{span}
\lim_{m\to\infty} {\rm span}(P(A,A^{-4n_1(m)},\ldots,A^{-4n_k(m)}))=\infty
\end{equation}
Let $\x=(x_1,\ldots,x_k)$, $\r=(r_1,\ldots,r_k)$, and $\x^{\r}=\prod_{i=1}^{k}x_i^{r_i}$.

\paragraph{Case 1}
Suppose $\displaystyle\lim_{m\to \infty} \nu(1,n_1(m),\ldots,n_k(m))=\infty$.
By \fullref{BLe}:
\begin{eqnarray*}
M(P(A,x_1,\ldots,x_k)) &=& \lim_{m\to\infty} M(P(A,A^{-4n_1(m)},\ldots,A^{-4n_k(m)}))\\
&=& \lim_{m\to\infty} M(P(A,(-A^{-4})^{n_{1}(m)},\ldots,(-A^{-4})^{n_{k}(m)})) \\
&=&\lim_{m\to\infty}M(\kb{L_{\n(m)}})=1
\end{eqnarray*}
using $M(\delta)=1$. \fullref{cpol} now implies:
$$P(A,x_1,\ldots,x_k)=A^{\alpha}\x^{\r_0}\prod_{j=1}^{\ell}\phi_j\left(A^{\alpha_j}\x^{\r_j}\right)$$
By \eqref{span}, some $\r_j\neq 0$.
Using \fullref{sixthroot}, evaluate $P(\xi,x_1,\ldots,x_k)$ at $\xi=e^{2\pi i/6}$:
$$\prod_{i=1}^kx_i\prod_{j=k+1}^N(-\xi^{-3})^{n_j}=
\xi^{\alpha}\x^{\r_0}\prod_{j=1}^{\ell}\phi_j\left(\xi^{\alpha_j}\x^{\r_j}\right)$$
This is a contradiction:   The right--hand side can be made zero for an
appropriate choice of nonzero $x_i$'s, whereas the left--hand side
will remain nonzero.

\paragraph{Case 2}
Suppose $\nu(1,n_1(m),\ldots,n_k(m))$ is bounded for all $m$.  By
\fullref{nu}, after passing to a subsequence, there exists
$\a=(a_0,\ldots,a_{k})\in \Z^{k+1}\setminus \{\mathbf{0}\}$ such that
$\a\cdot (1,n_1(m),\ldots,n_k(m))=0$. Without loss of generality, let
$a_{k} > 0$. We define:
$$Q(A,x_1,\ldots,x_{k-1})=P(A^{a_k},x_1^{a_k},\ldots,
x_{k-1}^{a_k},A^{4a_0}\prod_{i=1}^{k-1}x_i^{-a_i})$$
$\a\cdot (1,n_1(m),\ldots,n_k(m))=0$ implies $a_kn_k(m)=-a_0-\sum_{i=1}^{k-1} a_i n_i(m)$.  Thus,
$$Q(A,A^{-4n_1(m)},\ldots,A^{-4n_{k-1}(m)})=P(A^{a_k},A^{-4a_kn_1(m)},\ldots,A^{-4a_kn_k(m)}).$$
By \eqref{span}, this implies
\begin{equation}\label{spanQ}
\lim_{m\to\infty} {\rm span}(Q(A,A^{-4n_1(m)},\ldots,A^{-4n_{k-1}(m)}))=\infty
\end{equation}
For any $a\in \Z\setminus 0$, $M(f(t))=M(f(t^a))$, so we have:
\begin{align*}
M(Q(A,A^{-4n_1(m)}\!,\ldots,A^{-4n_{k-1}(m)})) &= M(P(A^{a_k},A^{-4a_kn_1(m)}\!,\ldots,A^{-4a_kn_k(m)})) \\
&=  M(P(A,A^{-4n_1(m)},\ldots,A^{-4n_k(m)})) 
\end{align*}

\paragraph{Case 2a}
Suppose $\displaystyle\lim_{m\to \infty}\nu\left(1,n_1(m),\ldots,n_{k-1}(m)\right)= \infty$.
By \fullref{BLe}, the previous equation implies:
\begin{eqnarray*}
M(Q(A,x_1,\ldots,x_{k-1}))&=&\lim_{m\to\infty}M(Q(A,A^{-4n_1(m)},\ldots,A^{-4n_{k-1}(m)}))\\
&=& \lim_{m\to\infty}M(P(A,A^{-4n_1(m)},\ldots,A^{-4n_k(m)}))\\
&=&\lim_{m\to\infty}M(\kb{L_{\n(m)}})=1
\end{eqnarray*}
\fullref{cpol} now implies:
$$Q(A,x_1,\ldots,x_{k-1}) = A^{\alpha}\x^{\r_0}\prod_{j=1}^{\ell}\phi_j\left(A^{\alpha_j}\x^{\r_j}\right)$$
By \eqref{spanQ}, some $\r_j\neq 0$.
Using \fullref{sixthroot}, evaluate $Q(\xi^{1/a_k},x_1,\ldots,x_{k-1})$ at $\xi=e^{2\pi i/6}$,
and we obtain the same contradiction as in Case 1.

\paragraph{Case 2b}
Suppose $\nu(1,n_1(m),\ldots,n_{k-1}(m))$ is bounded
for all $m$.  As we did at the start of Case 2, we reduce
$Q(A,x_1,\ldots,x_{k-1})$ to obtain a $(k-1)$-variable polynomial.  If
$\nu(1,n_1(m),\ldots,n_{k-2}(m))$ is bounded for all $m$, we again
reduce to obtain a $(k-2)$-variable polynomial.  Proceeding in a
similar manner, we may finally reach $(1,n_1(m))$ where $n_1(m)$ is
unbounded, so $\lim_{m\to \infty}\nu(1,n_1(m))=\infty$.  Whenever we
find $j$ such that
$\lim_{m\to\infty}\nu\left(1,n_1(m),\ldots,n_{k-j}(m)\right)= \infty$,
we proceed as in Case 2a to obtain a contradiction.  \end{proof}

\begin{example}\rm
This example demonstrates the reduction method in Case 2 of the proof
of \fullref{mainthm}.  For some $P(A,x_1,x_2,x_3)$ in our context, suppose:
$$\n(m)=(2m,3m^2+1,5m-2)$$
$$M(\kb{L_{\n(m)}})=M(P(A,(-A^{-4})^{2m},(-A^{-4})^{3m^2+1},(-A^{-4})^{5m-2}))$$
Since $(4,-5,0,2)\cdot (1,2m,3m^2+1,5m-2)=0$, then
$\nu(1,2m,3m^2+1,5m-2)$ is bounded for all $m$.
$$Q(A,x_1,x_2)=P(A^2,x_1^2,x_2^2,A^{16}x_1^{5})$$
$$Q(A,(A^{-4})^{2m},(A^{-4})^{3m^2+1}) = P(A^2,(A^{-8})^{2m},(A^{-8})^{3m^2+1},(A^{-8})^{5m-2})$$
Since $\displaystyle\lim_{m\to\infty}\nu(1,2m,3m^2+1)= \infty$:
$$M(Q(A,x_1,x_2))=\lim_{m\to\infty}M(Q(A,(A^{-4})^{2m},(A^{-4})^{3m^2+1}))=\lim_{m\to\infty}M(\kb{L_m})=1$$
\end{example}

\bibliographystyle{gtart}
\bibliography{link}

\begin{thebibliography}{}
\providecommand\bibmarginpar{\leavevmode\marginpar}
\def\urlstyle#1{{\tt #1}}

\bibitem{Boyd81}
\textbf{D\,W Boyd}, \href{http://dx.doi.org/10.1016/0022-314X(81)90033-0}
  {\emph{Kronecker's theorem and {L}ehmer's problem for polynomials in several
  variables}}, J. Number Theory 13 (1981) 116--121 \xox{MR}{602452}

\bibitem{ck05}
\textbf{A Champanerkar}, \textbf{I Kofman},
  \href{http://dx.doi.org/10.2140/agt.2005.5.1} {\emph{On the {M}ahler measure
  of {J}ones polynomials under twisting}}, Algebr. Geom. Topol. 5 (2005) 1--22
  \xox{MR}{2135542}

\bibitem{ckp}
\textbf{A Champanerkar}, \textbf{I Kofman}, \textbf{E Patterson},
  \href{http://dx.doi.org/10.1142/S021821650400355X} {\emph{The next simplest
  hyperbolic knots}}, J. Knot Theory Ramifications 13 (2004) 965--987
  \xox{MR}{2101238}

\bibitem{Kanenobu2}
\textbf{T Kanenobu}, \href{http://dx.doi.org/10.1007/BF01459137}
  {\emph{Examples on polynomial invariants of knots and links}}, Math. Ann. 275
  (1986) 555--572 \xox{MR}{859330}

\bibitem{Kanenobu1}
\textbf{T Kanenobu}, \emph{Infinitely many knots with the same polynomial
  invariant}, Proc. Amer. Math. Soc. 97 (1986) 158--162 \xox{MR}{831406}

\bibitem{Lackenby}
\textbf{M Lackenby}, \href{http://dx.doi.org/10.1112/S0024611503014291}
  {\emph{The volume of hyperbolic alternating link complements}}, Proc. London
  Math. Soc. $(3)$ 88 (2004) 204--224 \xox{MR}{2018964}\ \ With an appendix by
  Ian Agol and Dylan Thurston

\bibitem{Lawton}
\textbf{W\,M Lawton}, \href{http://dx.doi.org/10.1016/0022-314X(83)90063-X}
  {\emph{A problem of {B}oyd concerning geometric means of polynomials}}, J.
  Number Theory 16 (1983) 356--362 \xox{MR}{707608}

\bibitem{Lickorish}
\textbf{W\,B\,R Lickorish}, \emph{An introduction to knot theory}, Graduate
  Texts in Mathematics 175, Springer, New York (1997) \xox{MR}{1472978}

\bibitem{Schinzel}
\textbf{A Schinzel}, \emph{The {M}ahler measure of polynomials}, from: ``Number
  theory and its applications (Ankara, 1996)'', Lecture Notes in Pure and Appl.
  Math. 204, Dekker, New York (1999)  171--183 \xox{MR}{1661667}

\bibitem{SSW}
\textbf{D Silver}, \textbf{A Stoimenow}, \textbf{S Williams},
  \href{http://dx.doi.org/10.2140/agt.2006.6.581} {\emph{Euclidean {M}ahler
  measure and twisted links}}, Algebr. Geom. Topol. 6 (2005) 581--602

\bibitem{Watson}
\textbf{L Watson}, \emph{Any tangle extends to non-mutant knots with the same
  {J}ones polynomial}, to appear in J. Knot Theory Ramifications

\end{thebibliography}

\end{document}